\documentclass[12pt]{article}
\usepackage{amsmath,amssymb,amscd,latexsym,mathrsfs}
\usepackage{color}

\usepackage[all]{xy}

\newcommand{\colim}{{\rm colim}}

\newcommand{\hocolim}{{\rm hocolim}}
\newcommand{\coLim}{\underrightarrow{\lim}}
\newcommand{\dom}{{\rm dom\,}}
\newcommand{\cod}{{\rm cod\,}}

\newcommand{\Cat}{{\rm Cat}}
\newcommand{\Set}{{\rm Set}}
\newcommand{\sSet}{{\rm sSet}}

\newcommand{\Grp}{{\rm Grp}}
\newcommand{\Ab}{{\rm Ab}}
\newcommand{\Ob}{{\rm Ob}}
\newcommand{\Mor}{{\rm Mor}}

\newcommand{\NN}{{\,\mathbb N}}

\newcommand{\mC}{{\,\mathscr C}}
\newcommand{\mD}{{\,\mathscr D}}

\newcommand{\mG}{{\cal G}}

\newcommand{\mA}{{\mathscr A}}

\newcommand{\cA}{{\mathcal A}}
\newcommand{\mB}{{\mathcal B}}

\newcommand{\mX}{{\mathcal X}}
\newcommand{\fF}{{\mathfrak F}}

\newtheorem{theorem}{\bf Theorem}[section]
\newtheorem{lemma}[theorem]{\bf Lemma}
\newtheorem{proposition}[theorem]{\bf Proposition}
\newtheorem{corollary}[theorem]{\bf Corolary}
\newtheorem{definition}{\sc Definition}[section]
\newtheorem{example}[definition]{\sc Example}
\newtheorem{remark}[definition]{\sc Remark}

\def\leq{\leqslant}
\def\geq{\geqslant}

\def\MYstar{\mathop{\star}\limits}

\begin{document}

\begin{center}
{\large \bf Homotopy cofinality for Non-Abelian homology \\
of group diagrams
}
\\
\medskip
Ahmet A. Husainov\\
\end{center}

\begin{abstract}
We prove that a homotopy cofinal functor between small categories induces a weak equivalence between homotopy colimits of pointed simplicial sets. This is used to prove that the non-Abelian homology of a group diagram is isomorphic to the homology of its inverse  image under a homotopy cofinal functor.
This also made it possible to establish that the non-Abelian homology of group diagrams are invariant under the left Kan extension along virtual discrete cofibrations. With the help of these results, we have constructed a non-Abelian Gabriel-Zisman homology theory for simplicial sets with coefficients in group diagrams.
Moreover, we have generalized this theory to presheaves over an arbitrary small category, which play the role of simplicial sets.
Sufficient conditions are found for the isomorphism of non-Abelian homology of presheaves over a small category with coefficients in group diagrams.
It is proved that the non-Abelian homology of the factorization category for the small category is isomorphic to the Baues-Wirsching homology.
A condition is found under which a functor between small categories induces an isomorphism of non-Abelian Baues-Wirsching homology.
\end{abstract}

2000 Mathematics Subject Classification 55U10, 18C15, 18G10, 55P10, 55R35

Keywords: pointed simplicial sets, homotopy groups, simplicial groups, group diagrams, non-Abelian homology, small category homology, Baues-Wirsching homology, factorization category homology, Gabriel-Zisman homology, simplicial replacement, classifying spaces, homotopy colimit, homotopy cofinality, Grothendieck construction.

\tableofcontents

\section{Introduction}

In this paper, we solve the following two problems in the theory of non-Abelian homology for group diagrams:

\begin{itemize}
\item Let $S: \mC\to \mD$ be a homotopy cofinal functor between small categories and $\mG: \mD\to \Grp$ be a group diagram.
Will the canonical homomorphisms of non-Abelian homology $\colim^{\mC}_n\mG S\to \colim^{\mD}_n\mG$ be isomorphisms for all $n\geq 0$?
\item Let $S: \mC\to \mD$ be a discrete Grothendieck cofibration between small categories and $\mG: \mC\to \Grp$ be a group diagram. Is it true that the canonical homomorphisms $\colim^{\mC}_n \mG \to \colim^{\mD}_n Lan^S\mG$ are isomorphisms for all $n\geq 0$?
\end{itemize}

 Section 3 is devoted to the theorem that a composition with a homotopy cofinal functor preserves the homotopy colimit of a diagram taking values in the category of pointed simplicial sets. This was previously known for diagrams of non-pointed simplicial sets \cite[\S 9.4]{dwy1984}.
  We prove this theorem with a lemma on the connection between pointed and non-pointed homotopy colimits \cite{far1995}, \cite{far2004}, \cite{hir2003}.

 Section 4 begins with the definition of non-Abelian homologies for group diagrams $\mG: \mC \to \Grp$. For example, if $\mC$ is equal to the group $\Gamma$, then these homologies turn into non-Abelian homologies of $\Gamma$-groups \cite{ina2022}.
  A description is given of the basic formula used to calculate and study these homologies. A theorem is proved that if a functor $S: \mC\to \mD$ is homotopy cofinal, then for any group diagram $\mG: \mD\to \Grp$ the homomorphisms $\colim^{\mC}_n\mG S \to \colim^{\mD}_n \mG$ are isomorphisms for all $n\geq 0$.
   A similar assertion for homology of small categories with coefficients in Abelian categories with exact coproducts was proved by Oberst \cite{obe1968}.
 
 We also prove a theorem stating that if the functor $S: \mC\to \mD$ is a virtual discrete cofibration, then, for the left Kan extension $Lan^S$ along $S$, the formula 
$\colim^{\mC}_n\mG \cong \colim^{\mD}_n Lan^S \mG$, for every group diagram $\mG$ on $\mC$ and all $n\geq 0$.
As an example, a category with a unique decomposition of a morphism into a composition of a monomorphism and an epimorphism is given. The embedding functor of a subcategory of monomorphisms into this category will be a virtual discrete cofibration.

Section 5 considers non-Abelian homology of simplicial sets with coefficients in group diagrams, which are constructed similarly to the Gabriel-Zisman homology with coefficients in diagrams of Abelian objects. As an example, the non-Abelian Baues-Wirsching homology \cite{bau1985} is considered as the Gabriel-Zisman homology for the category nerve and it is proved that they are isomorphic to the non-Abelian homology of a contravariant natural system (Corollary \ref{corfact}). This was previously known for the Abelian Baues-Wirsching cohomology \cite[Theorem 4.4]{bau1985}. A large number of applications of Abelian cohomology and Abelian Gabriel-Zisman homology are given in \cite{gal2021}.
 
 We show for non-Abelian homology that instead of the category $\Delta$ one can substitute an arbitrary small category $\mD$ and instead of simplicial sets consider $\mD$-sets - presheaves over $\mD$. Many homology properties of simplicial sets hold for $\mD$-sets. Among these properties are some isomorphism conditions for non-Abelian homology for $\mD$-sets. For Abelian homology, this was done in \cite{X2022}.
 
 \section{Preliminaries}
 
Let us write out the main definitions and notation. The rest will be given in the course of the presentation.

\begin{itemize}
\item
$\mA^{\mC}$ is the category of functors from the small category $\mC$ to an arbitrary category $\mA$.
Functors $F$ from a small category $\mC$ to an arbitrary $\mA$ are called object diagrams of the category $\mA$ over $\mC$ and can be denoted as the family $\{F(c)\}_{c\in \mC}$.
\item $\Delta_{\mC}A$ is a functor $\mC\to \mA$ that takes constant values $A\in \Ob\mA$ on objects and $1_A$ on morphisms.
\item $\NN$ - set of non-negative integers.
\item $\Cat$ - category of small categories and functors.
\item
$\Set$ - category of sets and mappings.
\item
${\Grp}$ is the category of groups and homomorphisms.
\item $0$ - a group consisting of one element.
\item $\Delta$ - the category of finite linearly ordered sets $[n]= \{0, 1, \cdots, n\}$, $n\geq 0$,
and non-decreasing mappings. The category $\Delta$ is generated by morphisms of the following form:
\begin{enumerate}
\item $\partial^i_n: [n-1]\to [n]$ (for $0\leq i\leq n$) - increasing mapping whose image does not contain $i$,
\item
$\sigma^i_n: [n+1]\to [n]$ (for $0\leq i\leq n$) is a non-decreasing surjection that takes the value $i$ twice,
\end{enumerate}
\item
$\colim^{\mC}_n: \Grp^{\mC}\to \Grp$ ($\forall n\geq 0$) are functors assigning to each group diagram the values of non-Abelian homology of this diagram (see \cite{ X2023}), $\colim^{\mC}$ denotes the colimits of group diagrams,
\item
$\coLim^{\mC}_n: \mA^{\mC}\to \mA$ ($\forall n\geq 0$) - left satellites of the colimit functor from the category of diagrams on $\mC$ taking values in the abelian category $\mA$ with exact coproducts \cite[Application 2]{gab1967}.
\end{itemize}

Let $\Phi: \mC\to \mD$ be a functor between categories.
For each $d\in \Ob\mD$ the left fibre \cite[Appendix 2, \S3.5]{gab1967} or the comma category 
$\Phi$ over $d$ \cite[\S2.6]{mac2004} is the category $\Phi\downarrow d$, whose objects are the pairs 
$(c\in \Ob\mC, \beta\in \mD(S(c), d))$, and the morphisms $(c, \beta)\to (c', \beta')$ are 
given by $\alpha\in \mC(c,c')$ satisfying $\beta'\circ S(\alpha)= \beta$. A morphism in 
$\Phi\downarrow d$ is defined by a triple consisting of two objects and a morphism between them.

For an arbitrary category $\mA$, the inverse image functor $\Phi^*: \mA^{\mD}\to \mA^{\mC}$ is defined, which assigns to each diagram $F\in \mA^{\mD}$ the composition $F\Phi= F\circ \Phi\in \mA^{\mC}$, and to each natural transformation $\eta: F\to F'$ the natural transformation $\eta\Phi: F\Phi \to F'\Phi$ defined by the formula $(\eta\Phi)_{c}= \eta_{\Phi(c)}$, for all $c\in \Ob\mC$.
If $\mA$ is a cocomplete category, then the functor $\Phi^*$ has a left adjoint functor $Lan^{\Phi}: \mA^{\mC}\to \mA^{\mD}$, which is called the left Kan extensions.
For an arbitrary diagram $F\in \mA^{\mC}$, the diagram $Lan^{\Phi}F$ can be considered as taking values on the objects $d\in \Ob\mD$ equal to $Lan^{\Phi}F (d)= \coLim^{\Phi\downarrow d}FQ_d$, where $Q_d: \Phi\downarrow d\to \mC$ maps each object $S(c)\to d$ to an object $c\in \mC $, with an obvious extension to the morphisms \cite[\S10.3, (10)]{mac2004}.

In the case when $\Phi: \mC\to \mD$ is a full embedding of the small category $\mC$ in the category $\mD$, we will denote the category $\Phi\downarrow d$ by $\mC\downarrow d$ .

For example, for the Yoneda embedding $h^{\mD}: \mD\to \Set^{\mD^{op}}$ and the presheaf $X\in \Set^{\mD^{op}}$, the category $ h^{\mD} \downarrow X$ is denoted by $\mD\downarrow X$.

Objects of the category $\Set^{\mD^{op}}$, functors $\mD^{op}\to \Set$, are called $\mD$-sets, and natural transformations between them are called morphisms between $\mD$ -sets.

  A simplicial set is a functor $X: \Delta^{op}\to \Set$. The values of this functor on morphisms generating the category $\Delta$ are denoted by $d^n_i= X(\partial^i_n)$, $s^n_i= X(\sigma^i_n)$, for all $0\leq i\leq n$.

$\Delta[n]$, $n\geq 0$, denotes the \textit{standard} simplicial set defined as the functor $\Delta(-, [n]): \Delta^{op}\to \Set$ .

A simplicial mapping $X\to Y$ between simplicial sets is a natural transformation. A simplicial mapping $X\to Y$ is called a cofibration if for each $n\geq 0$ the mapping $X_n\to Y_n$ is an injection.

The category of simplicial sets is denoted by $\Set^{\Delta^{op}}$ or $\sSet$.

A simplicial set $X$ together with a distinguished point $x\in X_0$ is called pointed.
We will consider this point $x$ as a morphism $\tilde{x}: \Delta[0]\to X$ corresponding to this point by Yoneda's lemma.
The category of pointed sets is equal to $\Delta[0]\downarrow \Set^{\Delta^{op}}$ and is shortly denoted by $\sSet_*$.

We adhere to the theory of homotopy colimits described in \cite{bou1972}, \cite{goe2009}.

Let $\mC$ be a small category. Consider two functors
$$
\mC \xleftarrow{\partial_0} (\Delta\downarrow \mC)^{op} \xrightarrow{Q^{op}_{\mC}}
\Delta^{op}.
$$
The functor $\partial_0$ assigns to each singular simplex $\sigma = (c_0\to c_1 \to \cdots \to c_n)$ in a nerve of the category $\mC$ its origin $c_0$.
It can be defined using a mapping that associates the functor $\sigma: [n]\to \mC$ with the object $\sigma(0)\in \mC$.
And the functor $Q_{\mC}: \Delta\downarrow \mC \to \Delta$ takes $\sigma: [n]\to \mC$ to $[n]\in \Delta$.

For an arbitrary category $\cA$, the functor of \textit{simplicial replacement} 
$\coprod: \cA^{\mC}\to \cA^{\Delta^{op}}$ is defined by the formula 
$\coprod F = Lan ^{Q^{op}_{\mC}} (F\circ \partial_0)$.
The meanings of the simplicial replacement are detailed in \cite[Proposition 2.2]{X2022}.

In particular, if $\cA= \sSet_*$, then for $\mX\in \sSet_*^{\mC}$ the simplicial replacement consists of pointed simplicial sets $\coprod_n \mX= \vee_{c_0\to \cdots \to c_n} \mX(c_0)\in \sSet_*$. It gives a bisimplicial set
$$
\vee_{c_0\to \cdots \to c_n}\mX(c_0)_m.
$$
The homotopy colimit of the diagram $\mX: \mC\to \sSet_*$ can be defined as the diagonal of the simplicial replacement
$\hocolim^{\mC}\mX:= diag\coprod\mX$ \cite[Lemma XII.5.2]{bou1972}. It consists of sets
$$
  (\hocolim^{\mC}\mX)_n= \vee_{c_0\to \cdots \to c_n}\mX(c_0)_n.
$$
The marked vertex is obtained by identifying marked vertices from $\mX(c_0)$, $\forall c_0\in \Ob\mC$, and is a common point of $\vee_{c_0}\mX(c_0)_0$.
In each dimension, the degenerations of this vertex are identified.

\section{Cofinal functors and homotopy colimits}

This section is about homotopy colimits of pointed simplicial sets and homotopy cofinal functors introduced in \cite{hir2003}.

\begin{definition}
Let $S: \mC\to \mD$ be a functor between small categories.
It is called homotopy cofinal if, for every object $d\in \mD$, the nerve of category $d\downarrow S$ is non-empty and contractible (i.e., the geometric realization for the nerve of the category $d\downarrow S$ has homotopy groups like a point). We will call a functor $S$ homotopy coinitial if $S^{op}$ is homotopy cofinal.
\end{definition}

In \cite[\S 9.4]{dwy1984} a formula was given showing that if a functor $S: \mC\to \mD$ is homotopy cofinal, then for any diagram of simplicial sets $\mX$ over $\mD$ the canonical mapping $\hocolim^{\mC}\mX S \to \hocolim^{\mD}\mX$ is a weak equivalence.
The proof of this formula is published in Hirshhorn's book \cite[Theorem 19.6.7]{hir2003}.

We prove a similar assertion for diagrams of pointed simplicial sets.

\subsection{Pointed and unpointed homotopy colimits}

Let us consider a statement due to Farjoun \cite{far1995} (see also \cite[Proposition 18.8.4]{hir2003}) about the connection between unpointed and pointed homotopy colimits in the category of pointed simplicial sets $\sSet_*:= \Delta[0] \downarrow \Set^{\Delta^{op}}$.
We formulate it in the following lemma in a form convenient for use in proving the confinal functor theorem for diagrams in $\sSet_*$.

Let $U: \Delta[0]\downarrow \Set^{\Delta^{op}} \to \Set^{\Delta^{op}}$ be an forgetful functor that associates each pointed simplicial set $\Delta[0]\to Y$ with an unpointed simplicial set $ Y$.
For an arbitrary diagram $\mX: \mD\to \Delta[0]\downarrow \Set^{\Delta^{op}}$ we denote by $U\mX:= U^{\mD}(X)$ the composition $ U\circ \mX: \mD\to \Set^{\Delta^{op}}$.

We gave the definition of a pointed homotopy colimit in the preliminary information. This definition is based on the statement \cite[Lemma XII.5.2]{bou1972}.

In \cite[Notation 18.1.1]{hir2003} two new types of colimit were introduced and called them $\hocolim^{\mC}_*$ and $\hocolim^{\mC}$. The notation for the second of these types is busy, so we will denote these colimits by $Hocolim^{\mC}_*$ and $Hocolim^{\mC}$, respectively.

Analysis of the text \cite[Notation 18.1.1]{hir2003} shows that these operators have common domains and codomains:
$$
Hocolim^{\mC}_*,~ Hocolim^{\mC}: \sSet^{\mC}_*\to \sSet.
$$
These operators decompose into (not necessarily equal) two compositions along the main diagonal of the square
$$
\xymatrix{
\sSet^{\mC}_*
\ar[d]_{U^{\mC}} 
\ar[rr]^{\hocolim} && \ar[d]^U \sSet_* \\
\sSet^{\mC} \ar[rr]^{\hocolim} &&  \sSet
}
$$
We arrive at the following formulas for the introduced operators

\begin{itemize}
\item $Hocolim^{\mC}_*\mX:= U(\hocolim^{\mC}\mX$),
\item $Hocolim^{\mC}\mX:= \hocolim^{\mC}(U^{\mC}(\mX))=
\hocolim^{\mC}(U\mX)$.
\end{itemize}

The following lemma is auxiliary and is proven using Hirschhorn's statement \cite[Prop. 18.8.4]{hir2003}.

\begin{lemma}\label{lcodecar}
\begin{enumerate}
\item
For every diagram $\mX: \mD \to \Delta[0]\downarrow \Set^{\Delta^{op}}$ there is a quadrangle of morphisms in the category of simplicial sets
\begin{equation}\label{codecar}
\xymatrix{
B\mD^{op}\ar[d] \ar[rr] && \hocolim^{\mD} U\mX\ar[d]\\
{\Delta[0]} \ar[rr] && U(\hocolim^{\mD}\mX)
}
\end{equation}
whose top row is a natural cofibration.

\item

Homotopy colimit contained in (\ref{codecar})
diagrams
\begin{equation} \label{loc4}
\Delta[0] \overset{e}{\leftarrow} B\mD^{op} \overset{f}{\to} \hocolim^{\mD}U\mX
\end{equation}
equals 
$U(\hocolim^{\mD}\mX)$.

\item
A pair consisting of the simplicial set $U(\hocolim^{\mD}\mX)$ and the morphism ${\Delta[0]} \to U(\hocolim^{\mD}\mX)$ in the bottom row of (\ref{codecar}) can be considered as a pointed simplicial set.
\end{enumerate}
\end{lemma}
{\sc Proof.}
Let's start with the first statement:
In work \cite[Prop. 18.8.4]{hir2003} natural cofibration constructed
$B\mD^{op}\to \hocolim^{\mD}U\mX$ and natural isomorphism
\begin{equation}\label{isohir}
(hocolim^{\mD}U\mX)/(B\mD^{op})\cong U(hocolim^{\mD}\mX).
\end{equation}
We get a quadrilateral (\ref{codecar}),
whose top row is the natural cofibration $f$.
\begin{equation}\label{codecar2}
\xymatrix{
B\mD^{op}\ar[d]^e \ar[rr]^f && \hocolim^{\mD} U\mX\ar[d]\\
{\Delta[0]} \ar[rr]  && U(\hocolim^{\mD}\mX)\cong (\hocolim^{\mD}
U\mX/B\mD^{op})
}
\end{equation}

Let's move on to the second statement:
The homotopy colimit of a pushout (\ref{loc4}) will be equal to the colimit of a diagram consisting of the following pair of parallel arrows:
$$
B\mD^{op}  \underset{e}
{\overset{f} {\rightrightarrows}}
\Delta[0]\coprod \hocolim^{\mD} U\mX
$$
and therefore to the cokernel denoted by $\rm coker(e,f) $.

  The colimit of the diagram (\ref{loc4}) will be isomorphic
$\hocolim^{\mD}U\mX/ (B\mD^{op})$ which is equal to $U(\hocolim^{\mD}\mX)$
  according to the formula (\ref{isohir}).

The third statement is obvious.
  
\hfill$\Box$



\subsection{Preservation of a pointed homotopy colimit}

In this subsection we prove a theorem that will help us answer the questions formulated in the introduction.
  
\begin{theorem}\label{cofpointed}
Let $S: \mC\to \mD$ be a homotopy cofinal functor between small categories.
Then for every diagram $\mX$ of pointed simplicial sets over $\mD$ the canonical map $\hocolim^{\mC} (\mX S) \to \hocolim^{\mD}\mX$ is a weak equivalence.
\end{theorem}
{\sc Proof.} By Lemma \ref{lcodecar}, the space $U(\hocolim^{\mD}\mX)$ will be a homotopy colimit of the following diagram
$$
\Delta[0] \longleftarrow B\mD^{op} \rightarrowtail \hocolim^{\mD}U\mX.
$$
Now consider the morphism
$
\hocolim^{\mC} U\mX S \to \hocolim^{\mD} U\mX$, which, due to the homotopy confinality condition of the functor $S$, will be a weak equivalence, according to \cite[Theorem 19.6.7] {hir2003}. Since the nerve $B(d\downarrow S)$ is contractible for every $d\in \Ob\mD$, then by Quillen's theorem \cite[Theorem A]{qui1973} the mapping $B(S):B\mC \to B \mD$ will be a weak equivalence. In this case, the nerve $B(S^{op}\downarrow d)$ will also be contractible, and hence the mapping $B(S^{op}):B\mC^{op} \to B\mD^{op} $ would be a weak equivalence.

The naturalness of the mapping $\hocolim^{\mC} U\mX S\to \hocolim^{\mD} U\mX$ to $\mX$ leads to the following commutative diagram
$$
\xymatrix{
 B\mC= \hocolim^{\mC}U\Delta_{\mC}(*)S \ar[d] \ar[r] & \hocolim^{\mC}
 U\mX S \ar[d]\\
B\mD= \hocolim^{\mD} U\Delta_{\mD}(*) \ar[r] & \hocolim^{\mD} U\mX  
}
$$
where $\Delta_{\mD}(*)$ is a diagram over $\mD$ taking constant values equal to $*=\Delta[0]$.
Let's add to this diagram on the left a commutative diagram of nerves
$$
\xymatrix{
 B\mC^{op} \ar[d] \ar[r] & B{\mC}\ar[d]\\
B\mD^{op} \ar[r] & B{\mD}  
}
$$
The rows of this diagram are isomorphisms, acting for each $n\in\NN$ as $(c_0\gets \ldots \gets c_n) \mapsto (c_0 \to \ldots \to c_n)$, and similarly for the bottom row.

We obtain a natural transformation of diagrams into the category of pointed simplicial sets, shown with the help of three vertical arrows, each of which is a weak equivalence
$$
\xymatrix{
{\Delta[0]}\ar[d] & B\mC^{op} \ar[l] \ar[d] \ar[r] & \hocolim^{\mC} U\mX S 
\ar[d]\\
{\Delta[0]} & B\mD^{op} \ar[l] \ar[r] & \hocolim^{\mD} U\mX
}
$$
Colimits of the top and bottom rows lead to a commutative diagram
$$
\xymatrix{
{\Delta[0]}\ar[d] \ar[r]& U(\hocolim^{\mC}\mX S)  \ar[d]  &
 \hocolim^{\mC} U\mX S \ar[l] \ar[d]\\
{\Delta[0]} \ar[r] & U(\hocolim^{\mD}\mX)   & \hocolim^{\mD} U\mX \ar[l]
}
$$
whose vertical morphisms are weak equivalences.
From the third statement of Lemma \ref{lcodecar} it follows that the square on the left can be considered as a morphism of pointed homotopy colimits $\hocolim^{\mC}\mX S \to \hocolim^{\mD}\mX$, and therefore which will also be weak equivalence.
\hfill$\Box$
\section{Non-Abelian homology of group diagrams}

Non-Abelian homology for group diagrams is defined in \cite{X2023} as the 
cotriple derived functors of the colimit functor $\colim^{\mC}: \Grp^{\mC}\to \Grp$.
Nevertheless, we will show that to solve a number of problems in this theory, it suffices to define non-Abelian homologies of group diagrams using a simplicial change.
The technique for solving these problems is based on the formula for homotopy colimits of classifying spaces \cite[Theorem 5.2]{X2023}.

\subsection{Homotopy and homology of group diagrams}

Let $\mG: \mC\to \Grp$ be a group diagram over a small category $\mC$.
A simplicial replacement of a group diagram $\mG$ is a simplicial group $C_*(\mC, \mG)$ consisting of free products
$$
C_n(\mC, \mG)= \MYstar_{c_0\to \cdots \to c_n}\mG(c_0), n\geq 0.
$$
Its boundary operators $d^i_n: C_n(\mC,\mG)\to C_{n-1}(\mC,\mG)$ and degeneration $s^i_n: C_n(\mC,\mG)\to C_ {n+1}(\mC,\mG)$, for $0\leq i\leq n$, act on the elements of factors of free products by the formulas
\begin{multline}\label{fsimprep1}
d^i_n(c_0\xrightarrow{\alpha_1}c_1 \to 
\cdots \to c_{n-1} \xrightarrow{\alpha_n} c_n, x)=\\
=
\begin{cases}
(c_1\xrightarrow{\alpha_2} c_2\to \cdots \to c_{n-1} 
\xrightarrow{\alpha_n}c_n, \mG(\alpha_1)(x)),
& \text{ if $i=0$},\\
(c_0\xrightarrow{\alpha_1}c_1\to \cdots \to c_{i-1} \xrightarrow{\alpha_{i+1}\alpha_i}
c_{i+1}\to \cdots \xrightarrow{\alpha_n} c_n, x), & \text{ if $i>0$},
\end{cases}
\end{multline}
\begin{multline}\label{fsimprep2}
s^i_n(x, c_0\xrightarrow{\alpha_1}c_1 \to \cdots \to c_{n-1} \xrightarrow{\alpha_n} c_n)=\\
(c_0\xrightarrow{\alpha_1}c_1 \to \cdots \to c_i \xrightarrow{id} c_i
 \to \cdots \to c_{n-1} \xrightarrow{\alpha_n} c_n, x),
\end{multline}
where $x$ denotes an arbitrary element of the group $\mG(c_0)$.

\begin{definition}
Let $\mG: \mC\to \Grp$ be a group diagram. Its non-Abelian homologies are the groups
$$
\colim^{\mC}_n \mG := \pi_n(C_*(\mC, \mG)), n\geq 0.
$$
\end{definition}

The group $\colim^{\mC}_0\mG$ is equal to the colimit of the diagram $\mG$, it can be non-commutative. The groups $\colim^{\mC}_n\mG$ are commutative for $n\geq 1$ as homotopy groups of the simplicial group.

\begin{remark}
The homology diagrams $\mC\to \Ab$ are defined similarly, but instead of the coproduct $\star$ the direct sum $\oplus$ is used and all colimits are calculated in $\Ab$.
\end{remark}

Computation of non-Abelian homology is carried out using the following assertion, proved in the preprint \cite[Theorem 5.2]{X2023}.

\begin{theorem}\label{main2}
For an arbitrary group diagram $\mG: \mC\to {\Grp}$ and for all $n\geq 0$ there are natural isomorphisms
\begin{equation}\label{fmain2}
\colim^{\mC}_n \mG \cong
\pi_{n+1}(\hocolim^{\mC}B\mG).
\end{equation}
\end{theorem}

This allows one to study the non-Abelian homology of group diagrams as homotopy colimits of classifying spaces, methods for calculating which are developed in \cite{bro1987}, \cite{bro1990}, \cite{mik2010}.

\subsection{Cofinality and non-Abelian homology for group diagrams}

We give the following answer to the first question posed in the introduction.

\begin{theorem}\label{homoliso}
Let $S: \mC\to \mD$ be a homotopically cofinal functor between small categories. Then the canonical homomorphisms
$$
\colim^{\mC}_n \mG S \to \colim^{\mD} _n \mG
$$
are isomorphisms for any group diagram $\mG: \mD\to \Grp$ and for all $n\geq 0$.
\end{theorem}
{\sc Proof.}
Let $\mG: \mD\to \Grp$ be a group diagram. Consider the diagram of pointed simplicial sets $B\mG= \{B\mG(d)\}_{d\in\mD}$.
By Theorem \ref{cofpointed} simplicial mapping
\begin{equation}\label{miscel1}
\hocolim^{\mC}B\mG S \to \hocolim^{\mD}B\mG
\end{equation}
  will be a weak equivalence.
Consider the following commutative square
$$
\xymatrix{
\pi_{n+1}(\hocolim^{\mC}B\mG S)\ar[d] \ar[r] & 
\pi_{n+1}(\hocolim^{\mD}B\mG)\ar[d]\\
\colim^{\mC}_n \mG S \ar[r] & \colim^{\mD}_n \mG  
}
$$
for $n\geq 0$.
The homomorphism in the top row is invertible because the simplicial mapping (\ref{miscel1}) is a weak equivalence. The columns of the square contain isomorphisms by Theorem \ref{main2}. Hence the bottom row is an isomorphism.
\hfill$\Box$

\subsection{Formula for group diagram homology}

In this subsection, we prove the auxiliary Proposition \ref{hoclan} and Theorem \ref{discvirt} on the second question posed in the introduction.

\begin{definition}
Let $\mB$ be a category. It is called finally discrete if each connected component of the category $\mB$ contains a final object.
\end{definition}

\begin{lemma}\label{almdis}
Let $\mB$ be a category. It will be finally discrete if and only if it contains a discrete reflective subcategory.
\end{lemma}
(We adhere to the definition of reflectivity from \cite[Definition 3.5.2]{bor1994}.)

\begin{lemma}\label{holandiscr}
Let $\mC$ be a small finally discrete category and let $fin(\mC)$ be a discrete subcategory in $\mC$ consisting of final objects from connected components.
Then for any diagram $F: \mC\to \sSet_*$ there exists a natural (in $F\in \sSet^{\mC}_*$) weak equivalence
$$
   \hocolim^{\mC}F \approx \bigvee_{a\in fin(\mC)} F(a) .
$$
\end{lemma}
{\sc Proof.} Consider an embedding of a discrete subcategory consisting of final objects $fin(\mC)\subseteq \mC$. For each $c\in \Ob\mC$ the category $c\downarrow fin(\mC)$ consists of a single object, and so $B(c\downarrow fin(\mC))$ is contractible.
By Theorem \ref{cofpointed} $\hocolim^{fin(\mC)}F|_{fin(\mC)} \approx \hocolim^{\mC}F$, and the pointed homotopy colimit in the discrete category $fin( \mC)$ will be equal to $colim^{fin(\mC)}F|_{fin(\mC)}$ and hence $\bigvee_{a\in fin(\mC)} F(a)$.
\hfill$\Box$

For arbitrary small categories and a functor $S: \mC\to \mD$ between them, for any diagram $F: \mC\to \mA$ into the cocomplete category $\mA$, there exists a left Kan extension \cite[\S 10.3]{mac2004}
$$
Lan^S F (d)= \colim^{S\downarrow d} FQ_d
$$
where $Q_d: S\downarrow d\to \mC$ is a forgetful functor, $(S(c)\to d) \mapsto c$.

\begin{definition}
A functor $S: \mC\to \mD$ is called a virtual discrete cofibration if for each $d\in \Ob\mD$ the category $S\downarrow d$ is finally discrete.
\end{definition}

In particular, every discrete Grothendieck precofibration satisfies this condition.
The definition of the Grothendieck precofibration can be taken from Quillen's article \cite[\S1, Page 93]{qui1973}.

The term ``cofibration'' comes from Lemma \ref{almdis}.
The adjective ``virtual'' is necessary because the discrete Grothendieck precofibration also requires that the subcategory $fin(S\downarrow d)$ be equal to $S^{-1}(d)$.

\begin{proposition}\label{hoclan}
Let $S: \mC\to \mD$ be a virtual discrete cofibration between small categories. Then for any diagram of pointed simplicial sets $F: \mC\to \sSet_*$ there exists a natural weak equivalence
$$
\hocolim^{\mC} F \approx \hocolim^{\mD} Lan^S F,
$$
moreover, on objects the diagram $Lan^S F: \mD\to \sSet_*$ takes the values $Lan^S F(d) = \bigvee\limits_{\beta\in fin(S\downarrow d)} FQ_d(\beta)$.
\end{proposition}
{\sc Proof.} This follows from the formula
$$
\hocolim^{\mC}F \approx \hocolim^{\mD}
\{\hocolim^{S\downarrow d}FQ_d\}_{d\in \mD}~,
$$
due to Dwyer and Kan \cite[\S 9.8]{dwy1984} and applying Lemma \ref{holandiscr} for the category $S\downarrow d$ and the diagram $FQ_d$ for all $d\in \mD$.
\hfill$\Box$

The following theorem is used to study the connections between non-Abelian homologies of group diagrams over small categories.

\begin{theorem}\label{discvirt}
Let $S: \mC\to \mD$ be a virtual discrete cofibration between small categories.
Then for an arbitrary group diagram $\mG: \mC\to \Grp$ there exists an isomorphism of non-Abelian homology
\begin{equation}\label{fdiscvirt}
\colim^{\mC}_n \mG \cong
\colim^{\mD}_n Lan^S \mG,
~\forall n\geq 0,
\end{equation}
where $Lan^S\mG(d)$ for each $d\in \Ob\mD$ equals $\MYstar\limits_{\beta\in fin(S\downarrow d)} \mG Q_d(\beta)\} $ and $Lan^S\mG$ assigns to each morphism $\beta: d\to d'$ of the category $\mD$ a homomorphism defined on the factors of the free product as
$$
(\beta\in fin(S\downarrow d), g\in \mG Q_d(\beta)) \mapsto
(\phi(\alpha\beta), FQ_{d'}(\phi_{\alpha\beta})(g)).
$$
Here $\phi(\gamma)$, for $\gamma\in S\downarrow d'$, denotes the final object of the connected component of the category $S\downarrow d'$ containing $\gamma$, and
$\phi_{\gamma}: \gamma\to \phi(\gamma)$ is the only morphism to this final object.
\end{theorem}
{\sc Proof.} We will prove it using Theorem \ref{main2}.
Consider the diagram of pointed simplicial sets $B\mG: \mC\to \sSet_*$. Proposition \ref{hoclan} gives a weak equivalence
$$
\hocolim^{\mC} B\mG \approx \hocolim^{\mD}
\{\bigvee\limits_{\beta\in fin(S\downarrow d)} B\mG Q_d(\beta)\}_{d\in \mD}~.
$$
By Whitehead's theorem (see, \cite[Theorem 4.0]{fie1984} and \cite[Lemma 5.1]{X2023}), there is a weak equivalence $\bigvee_{j\in J}BG_j \approx B(\MYstar_{j\in J}G_j)$, for an arbitrary family of groups $(G_j)_{j\in J}$. We arrive at a weak equivalence
$$
\hocolim^{\mC} B\mG \approx \hocolim^{\mD}
\{B(\MYstar\limits_{\beta\in fin(S\downarrow d)} \mG Q_d(\beta))\}_{d\in \mD}~.
$$
Theorem \ref{main2}, after applying the functor $\pi_{n+1}$ to both parts of this weak equivalence, leads to the formula to be proved.
\hfill$\Box$

\begin{example}
Let $\mD$ be a small category in which each morphism admits a unique decomposition into the composition of an epimorphism and a monomorphism.
For example, $\mD$ is equal to the category of simplices $\Delta$ or the category of cubes $\Box$ \cite{X2022}. Denote by $\mD_+\subseteq \mD$ the subcategory consisting of monomorphisms with the set of objects $\Ob\mD_+= \Ob\mD$.
Denote by $J: \mD_+ \to \mD$ the embedding functor.
Just as \cite[Lemma 7.1]{X2022}, one can prove that for every object $d\in \mD$ every connected component of the category $d\downarrow J$ has an initial object, and hence the functor $J^{op} : \mD^{op}_+ \to \mD^{op}$ is a discrete virtual cofibration.
The set of final objects of the category $J^{op}\downarrow d$ will consist of epimorphisms $d\twoheadrightarrow d'$.
According to Theorem \ref{discvirt} for any diagram $\mG: \mD^{op}_+ \to \Grp$ the formula (\ref{fdiscvirt}) leads to an isomorphism:
$$
\colim^{\mD^{op}_+}_n \mG \cong
\colim^{\mD^{op}}_n
\{\MYstar\limits_{d\twoheadrightarrow d'} \mG(d')\}_{d\in \mD}~,
~\forall n\geq 0.
$$
\end{example} 

\section{Non-Abelian Gabriel-Zisman homology for $\mD$-sets}
 
 Let $\mD$ be a small category.
  A $\mD$-set or a presheaf of sets over $\mD$ is a functor $X: \mD^{op}\to \Set$. Morphisms between $\mD$-sets are natural transformations.
  We introduce non-Abelian homology of $\mD$-sets with coefficients in group diagrams, defined similarly to the Abelian homology of Gabriel and Zisman for simplicial sets \cite[Appendix 2]{gab1967}.
  We prove that for every $\mD$-set morphism $f: X\to Y$, the Kan extension of the group diagram on the $\mD$-set $X$ along $(\mD\downarrow f)^{op}$ preserves non-Abelian Gabriel-Zisman homology and we find conditions under which the composition of a group diagram on $Y$ with the functor $(\mD\downarrow f)^{op}$ does not change these homology.
 
 \subsection{Non-Abelian homology of $\mD$-sets with coefficients}
 
 Gabriel and Zisman \cite[Appendix 2, \S4]{gab1967} introduced homology $H_n(X, L)$ of a simplicial set $X$ with coefficients in the diagram $L:(\Delta\downarrow X)^{op} \to \mA$ of objects in the Abelian category $\mA$ with exact coproducts.
   And it was proved that these homologies are isomorphic to the values on $L$ of the left satellites of the colimit functor $\coLim^{(\Delta\downarrow X)^{op}}_n: \mA^{(\Delta\downarrow X)^{op }}\to \mA$ \cite[Appendix 2, Proposition 4.2]{gab1967}.
  The homology of $H_n(X, L )$ was defined as the homology of the simplicial object $Lan^{Q^{op}_X}L \in \mA^{\Delta^{op}}$ constructed as a left Kan extension along a functor, dual to the forgetful functor $Q_X: \Delta\downarrow X\to \Delta$ defined on objects as $Q_X(\Delta[n]\to X) \mapsto [n]$.
We prove the assertion from \cite[Appendix 2, Proposition 4.2]{gab1967} for the case when the role of the category $\Delta$ is played by an arbitrary small category $\mD$, and instead of the abelian category $\mA$ the category of groups $\Grp$.

 \begin{definition}
Let $\mD$ be a small category, $\Set^{\mD^{op}}$ be the category of $\mD$-sets. The non-Abelian homology groups of the $\mD$-set $X$ with coefficients in the diagram $\mG: (\mD\downarrow X)^{op}\to \Grp$ are defined by the formula
$$
H_n(X, \mG):= \colim^{\mD^{op}}_n Lan^{Q^{op}_X}\mG,
\text{for $n\geq 0$}.
$$
\end{definition}
These groups are commutative for $n\geq 1$.

Denote by $C^{\mD}_*(X,\mG): \mD^{op}\to \Grp$ the presheaf of groups over $\mD$ equal to $Lan^{Q^{op}_X}\mG$. The following assertion was proved for diagrams of Abelian groups in \cite[Lemma 14]{X2019}, and for diagrams of objects of the Abelian category with exact coproducts in \cite{X2022}. We need it for the case of group diagrams, which may not be commutative.
 
\begin{proposition}\label{contralan}
Let $\mG: (\mD\downarrow X)^{op}\to \Grp$ be a group diagram.
  Then there are natural isomorphisms
$$
\colim^{(\mD\downarrow X)^{op}}_n \mG \cong H_n(X, \mG)
$$
for all $n\geq 0$.
The functor $C^{\mD}_*(X, \mG): \mD^{op}\to \Grp$ associates with each $a \in \mD$ the free product $C^{\mD}_a(X, \mG)= \MYstar_{z\in X(a)}\mG(z)$, and each morphism $\alpha: b\to a$ has a homomorphism
$$
C^{\mD}_*(X, \mG) (\alpha): \MYstar\limits_{z\in X(a)} \mG(z) \to
  \MYstar\limits_{z\in X(b)} \mG(z),
$$
acting on elements of free factors as
\begin{multline*}
(x\in X(a), g\in \mG(x)) \mapsto\\
(X(\alpha)(x)\in X(b), \mG(x\xrightarrow{\alpha} X(\alpha)(x))(g)\in \mG(X(\alpha)(x))).
\end{multline*}
\end{proposition}
{\sc Proof} is given in \cite[Proposition 3.2]{X2022}, and differs in that in our assertion the category of groups is involved instead of the Abelian category, which means that the coproduct turns from a direct sum into a free product. The formula for $(Lan^{Q^{op}_X}\mG)(\alpha)$ is shown in \cite[Example 3.1]{X2022}.
\hfill$\Box$

\begin{corollary}\label{gzisos}
Let $\mG: (\Delta\downarrow X)^{op}\to \Grp$ be a group diagram. Then
$$
\colim^{(\Delta\downarrow X)^{op}}_n \mG \cong H_n(X, \mG)
\cong \pi_n(Lan^{Q^{op}_X}\mG).
$$
\end{corollary}
{\sc Proof.} The first isomorphism follows from Proposition \ref{contralan}, and the second from the formula $\colim^{\Delta^{op}}G\cong \pi_n(G)$, which is valid for an arbitrary simplicial group $G$ \cite[Corollary 5.3]{X2023}.
\hfill$\Box$
\begin{remark}
The first corollary isomorphism \ref{gzisos} is a non-Abelian version of \cite[Appendix 2, Proposition 4.2]{gab1967}, and Proposition \ref{contralan} generalizes this isomorphism to an arbitrary small category $\mD$ instead of $\Delta$.
\end{remark}

\subsection{Homology of direct and inverse images}

Let $f: X\to Y$ be a morphism of $\mD$-sets. For any contravariant group system $\mG: (\mD\downarrow X)^{op}\to \Grp$ on $X$ there exists a Kan extension $Lan^{(\mD\downarrow f)^{op}}\mG : (\mD\downarrow Y)^{op}\to \Grp$ on $Y$, which is called the direct image of the contravariant system $\mG$.
And for any contravariant system of groups $\mG$ on $Y$ there exists a contravariant system on $X$ equal to its composition $\mG\circ (\mD\downarrow f)^{op}$ with the functor $(\mD\downarrow f)^{op}$. It is called inverse image.

We will prove that Kan's extension along $(\mD\downarrow f)^{op}$ does not change the non-Abelian homology.
And in the case of composition with $(\mD\downarrow f)^{op}$, there are some sufficient conditions for these homology to not change.

\begin{lemma}
Let $f: X\to Y$ be a morphism of $\mD$-sets. Then the functor $(\mD\downarrow f)^{op}: (\mD\downarrow X)^{op} \to (\mD\downarrow Y)^{op}$ is a discrete cobundle.
\end{lemma}
{\sc Proof.}
In \cite[Lemma 3.4]{X2022} it is proved that for each $\tilde{y}\in \Ob(\mD\downarrow Y)$ each connected component of the category $\tilde{y}\downarrow (\mD\downarrow f)$ has an initial object. It follows that the category $(\mD\downarrow f))^{op}\downarrow \tilde{y}\cong (\tilde{y}\downarrow (\mD\downarrow f))^{op}$ is finally discrete , and hence the functor $(\mD\downarrow f)^{op}$ is a discrete cobundle.
\hfill$\Box$

\begin{corollary}\label{dliso}
For any $\mD$-set morphism $f: X\to Y$ and any contravariant group system $\mG: (\mD\downarrow X)^{op} \to \Grp$ there exists a natural non-abelian homology isomorphism
$$
H_n(X, \mG) \xrightarrow{\cong} H_n(Y, Lan^{(\mD\downarrow f)^{op}}\mG), \text{for all $n\geq 0$.}
$$
  \end{corollary}
{\sc Proof.} This follows from Theorem \ref{discvirt} and Proposition \ref{contralan}.
\hfill$\Box$

Let $\overleftarrow{f}({y})$ be a $\mD$-set defined by a Cartesian square in the category of $\mD$-sets
\begin{equation}\label{invimage}
\xymatrix{
X \ar[r]^f & Y \\
\overleftarrow{f}({y}) \ar[r] \ar[u]_{f_y}
& h^{\mD}_d\ar[u]_{\widetilde{y}}
}
\end{equation}
This $\mD$-set $\overleftarrow{f}(y)$ is called the {\it inverse fibre} over the element $y\in Y(d)$.
It is easy to see that there is an isomorphism of categories
$(\mD\downarrow f)\downarrow \widetilde{y}\cong \mD\downarrow \overleftarrow{f}(y)$.
Hence, if the nerve of the category $\mD \downarrow \overleftarrow{f}(y)$ is weakly equivalent to a point, then the functor $(\mD\downarrow f)^{op}$ is homotopy cofinal. From Theorem \ref{homoliso} we obtain the following sufficient condition for the isomorphism of non-Abelian homology of $\mD$-sets:

\begin{corollary}\label{dhiso}
Let $f: X\to Y$ be a morphism of $\mD$-sets.
If the nerves of the categories $\mD\downarrow \overleftarrow{f}(y)$ are contractible for all $y\in Y(d)$, then for any group diagram $\mG: (\mD\downarrow Y)^{op} \to \Grp$ canonical homomorphisms $H_n(X, \mG\circ (\mD\downarrow f)^{op}) \to
H_n(Y, \mG))$ are isomorphisms for all $n\geq 0$.
\end{corollary}

\begin{remark}
Our technique allows us to construct non-Abelian Andre homology
for $\mD$-sets. They differ from the Gabriel-Zisman homology in that instead of the contravariant group system $\mG: (\mD\downarrow X)^{op}\to \Grp$, the group diagram $\mG: \mD \to \Grp$ is taken.
Thus, for defining the non-Abelian Andre homology for the pair of categories $\Set^{\mD^{op}}\supset \mD$ and the functor $\mG: \mD \to \Grp$ considered in the first part of the monograph \cite{and1967}, we can take the formula
$$
H^A_n(X, \mG):= \colim^{\mD\downarrow X}_n\mG Q_X.
$$
The homology groups $H^A_n(h_d, \mG)$ for $\mD$-sets $h_d$, $d\in \Ob\mD$, will be equal to $0$ for all $n>0$, and equal to $\mG({h_d})$ for $n=0$.
For a group diagram $\mG$ on $\mD$ consisting of isomorphisms, the non-Abelian Andre homology $H^A_n(X, \mG)$ will be isomorphic to the Gabriel-Zisman homology $H_n(X, (\mG Q_X)^{- 1})$ with coefficients in the group diagram $(\mG Q_X)^{-1}: (\mD\downarrow X)^{op}\to \Grp$ obtained by homomorphism inversion.
\end{remark}

\section{Non-Abelian Baues-Wirsching homology}

This section is devoted to non-Abelian homology of contravariant natural systems of groups introduced similarly to the cohomology of small categories with coefficients in natural systems of Abelian groups \cite{bau1985}.

\subsection{Category of factorizations}

Let $\mC$ be a small category. The factorization category $\fF\mC$ \cite{bau1985} is defined as follows.
Its objects are the morphisms $f\in \Mor\mC$. Its morphisms $f\to g$ are given by the pairs $(\alpha,\beta)\in \Mor\mC^{op}\times \Mor\mC$ making the diagram commutative in the category $\mC$:
\begin{equation}\label{sqmor}
\xymatrix{
\circ \ar[r]^{\beta}& \circ \\
\circ \ar[u]^f & \circ \ar[l]^{\alpha} \ar[u]_g
}
\end{equation}
We will denote morphisms by $f \xrightarrow{(\alpha, \beta)}g$.
The identity morphism of the object $f\in \Ob\fF\mC$ is equal to $f \xrightarrow{(id_{\dom f}, id_{\cod f})} f$.
The morphism composition $f \xrightarrow{(\alpha, \beta)}g \xrightarrow{(\alpha', \beta')}h$ is defined as $f \xrightarrow{(\alpha\alpha', \beta' \beta )}h$.

We note that $\fF(\mC^{op})\cong \fF\mC$, where the isomorphism replaces all morphisms in the commutative diagram (\ref{sqmor}) with dual ones.

\begin{proposition}\label{wefrac}
The functor $\cod: \fF\mC \to \mC$ is homotopy coinitial, and the functor $\dom: (\fF\mC)^{op}\to \mC$ is homotopy cofinal.
In particular, the nerves of these categories are weakly equivalent.
\end{proposition}
{\sc Proof.} Let us prove homotopy coinitiality for the functor $\cod: \fF\mC \to \mC$. After that, substituting the category $\mC^{op}$ instead of $\mC$, we obtain the homotopy cofinality of the functor $\dom$.
 
Consider the category $\cod \downarrow c$, for $c\in \mC$. For each of its objects $(\alpha, f)$ there is a morphism $(1_{\dom\alpha},f)$ into the object $(f\circ\alpha, 1_c)$ shown in the diagram
$$
\xymatrix{
\dom\alpha \ar[rr]^{\alpha} && \cod\alpha \ar[d]_f  \ar[r]^f & c\\
\dom\alpha= \dom f\alpha \ar[u]_{1_{\dom\alpha}} \ar[rr]_{f\circ\alpha} && c \ar[ru]_{1_c}
} 
$$
The objects $(\beta, 1_c)$ constitute a full subcategory in $\cod\downarrow c$ equal to the Cartesian product $(\mC\downarrow c)^{op}\times\{1_c\}$. The morphism $(1_{\dom\alpha}, f)$ is the universal arrow from $(\alpha, f)\in \Ob(\cod\downarrow c)$ to the full embedding functor $in_c: (\mC\downarrow c )^{op}\times\{1_c\}\subseteq \cod\downarrow c$ in the sense of \cite[\S 3.1]{mac2004}, that is, for every morphism $(\alpha,f)\to (\beta , 1_c)$ there is a unique morphism $(f\alpha, 1_c)\to (\beta, 1_c)$ making the diagram commutative
$$
\xymatrix{
(\alpha, f) \ar[rd]_{\forall (g,f)} \ar[rr]^{(1_{\dom\alpha}, f)} && (f\alpha, 1_c)
\ar[ld]^{\exists!(g,1_c)}\\
& (\beta, 1_c)
}
$$
This implies \cite[\S4.1, Theorem 2]{mac2004} that the functor $i_c$ has a left adjoint fuector $l_c$, and hence the nerves of the categories $cod\downarrow c$ and $(\mC\downarrow c)^ {op}$ are homotopy equivalent. The category $\mC\downarrow c$ is homotopy equivalent to a point. Hence the functor $\cod$ is homotopy
coinitial.

\hfill$\Box$

To each functor $S: \mC\to \mD$ there corresponds a functor
$\fF S: \fF\mC \to \fF\mD$ defined on objects as
$\alpha\in \Ob\fF\mC \mapsto S(\beta)\in \Ob\fF\mD$, and on morphisms as
$$
(\alpha \xrightarrow{(u,v)} \alpha') \mapsto
  (S\alpha \xrightarrow{(Su,Sv)} S\alpha').
$$

To study the comma category $\fF S \downarrow \alpha$, $\alpha\in \Ob\fF\mD$, we introduce the category $S\langle \alpha\rangle$ whose objects are the pairs of morphisms $(u,v )$ from $\mD$ equal to the decompositions of the morphism $\alpha\in \Mor\mD$ into the composition $\dom \alpha \xrightarrow{u} S(d) \xrightarrow{v} \cod\alpha$.
The morphisms of $(u,v)\xrightarrow{\beta} (u', v')$ in the category $S\langle\alpha\rangle$ are the commutative diagrams
$$
\xymatrix{
\dom \alpha \ar[dr]_{u'} \ar[r]^u & S d \ar[d]^{S\beta} \ar[r]^v & \cod \alpha\\
& S d' \ar[ur]_{v'}
}
$$

The next assertion shows that the left fibers of the functor $\fF S$ are categories of factorizations.
\begin{proposition}\label{factfibres}
Let $S: \mC\to \mD$ be a functor between small categories.
Then for each $\alpha \in \Ob\fF\mD$ the category $(\fF S) \downarrow \alpha$ is isomorphic to $\fF(S\langle\alpha\rangle).$
\end{proposition}
{\sc Proof.} We construct a functor that assigns to each object $(S\beta \xrightarrow{(u,v)}\alpha) \in \Ob(\fF S\downarrow\alpha)$ an object $(u, vS\beta) \xrightarrow{\beta} (S(\beta)u, v)$ of the category $\fF(S\langle\alpha\rangle)$.
The action of a functor on morphisms $(u,v) \xrightarrow{(\gamma, \gamma')} (u',v')$ of the category $\fF S\downarrow \alpha$ can be shown using the following commutative diagrams
$$
\xymatrix{
S(\beta)\ar[d]_{(\gamma,\gamma')} \ar[r]^{(u,v)} & \alpha\\
S(\beta') \ar[ru]_{(u',v')}
}
\qquad
\xymatrix{
\ar@{}[d]_{\mapsto} \\
&
}
\xymatrix{
(u, vS(\beta)) \ar[r]^{\beta} & (S(\beta)u, v) \ar[d]^{\gamma'} \\
(u', v'S(\beta')) \ar[u]_{\gamma} \ar[r]_{\beta'} & (S(\beta')u', v') 
}
$$
The diagram on the left shows a morphism in the category $\fF S\downarrow \alpha$. On the right is the corresponding morphism from $\fF(S\langle\alpha\rangle)$.
If we add the morphism $(u',v')\xrightarrow{(\zeta, \zeta')} (u'',v'')$ on the left, then on the right we get one more rectangle adjoining from below.
On the right we get two commutative rectangles, proving the commutativity of the enclosing rectangle.
Thus, the composition of morphisms of the category $\fF S\downarrow \alpha$ goes over to composition, and hence the constructed mapping is functorial. It is clear that this mapping is bijective.
Therefore, it will be an isomorphism.
\hfill$\Box$

\subsection{Non-Abelian homology of the factorization category}

Baues and Wirsching \cite[Definition 1.4]{bau1985} introduced the cohomology of the small category with coefficients in the natural system of Abelian groups as the cohomology of the nerve of this category with coefficients in the diagram, which is equal to the composition 
$D\circ\delta$ of the natural system $D: \fF\mC \to \Ab$ and the functor $\delta: \Delta\downarrow B\mC \to \fF\mC$.
We define the non-Abelian Baues-Wirsching homology as the Gabriel-Zisman homology of the nerve $\mB\mC$ with coefficients in the group diagram $\mG\circ\delta^{op}: (\Delta\downarrow\mB\mC)^{ op} \to \Grp$.

Instead of the category of singular simplices of the nerve $\Delta\downarrow B\mC$, consider the isomorphic category $\Delta\downarrow \mC$ of functors defined on finite linearly ordered sets with values in the category $\mC$.
The functor $\delta: \Delta\downarrow \mC\to \fF\mC$ assigning to each functor $\sigma: [m]\to \mC$ an object $\sigma(0\leq m)$ of the category $\fF\mC $.
And it maps to each morphism
$$
\xymatrix{
[m]\ar[rd]_{\alpha} \ar[rr]^{\sigma} && \mC\\
&[n]\ar[ru]_{\tau}
}
$$
the morphism $\sigma(0\leq m)\xrightarrow{(\tau(0\leq \alpha(0)), \tau(\alpha(m)\leq n))}\tau(0\leq n)$ in the category $\fF\mC$.

\begin{lemma}\cite[Lemma 7.1]{X20232}\label{confBW}
For each $f\in \Ob\fF\mC$, the nerve of $\delta\downarrow f$ is weakly equivalent to a point.
\end{lemma}
{\sc Proof.} The idea of the proof is as follows.
Consider the category $S\langle\alpha\rangle$ for the functor $S=Id:\mC \to \mC$. The category $Id\langle\alpha\rangle$ will have initial and final objects, and hence its nerve is contractible.
Let us now prove that the category $\delta\downarrow f$ is isomorphic to the category $\Delta\downarrow Id\langle f\rangle$.

To construct the required isomorphism $\delta\downarrow f \cong \Delta\downarrow Id\langle f\rangle$, we note that each object of the category $\delta\downarrow f$ is given by the functor $[m]\xrightarrow{\sigma}\mC$ and the morphism $\delta(\sigma)\xrightarrow{(u,v)} f$ of the category $\fF\mC$, which is described by the following commutative diagram in $\mC$:
$$
\xymatrix{
\sigma(m) \ar[r]^v & \cod f\\
\sigma(0)\ar[u]^{\sigma(0\leq m)} & \dom f \ar[l]_u \ar[u]_f
}
$$
Associate this object with the functor $[m]\to Id\langle f\rangle$ which takes each element of $i\in [m]$ into an object of the category $Id\langle f\rangle$ equal to the decomposition of $f$ into composition
$$
\dom f \xrightarrow{g_i} \sigma(i) \xrightarrow{h_i} \cod f,
$$
where $g_i$ equals the composition $\dom f \xrightarrow{u} \sigma(0) \xrightarrow{\sigma(0\leq i)} \sigma(i)$ and $h_i$ is the composition $\sigma(i ) \xrightarrow{\sigma(i\leq m)} \sigma(m) \xrightarrow{v} \cod f$.
This comparison leads to a bijection between objects of categories $\delta\downarrow f$ and $Id\langle f\rangle$. Morphisms between objects of the category $\delta\downarrow f$ are given by morphisms from the category $\Delta$. The same is true for morphisms of the category $\Delta\downarrow Id\langle f\rangle$. This allows the bijection between objects to be extended to an isomorphism of categories.
The nerve of the category $Id\langle f\rangle$ is weakly equivalent to a point, and hence the nerve of the category $\Delta\downarrow Id\langle f\rangle$ is contractible.
Hence the nerve of the category $\delta\downarrow f$ is contractible.

\hfill$\Box$

A contravariant natural system of groups on a small category $\mC$ is a functor 
$\mG: (\fF\mC)^{op}\to \Grp$.
Following Baues and Wirsching \cite[Definition 1.4]{bau1985}, for the dual case we give the following 
\begin{definition}
 The non-Abelian homology groups $H^{BW}_n(\mC, \mG)$ of the category $\mC$ with coefficients in the contravariant natural group system $\mG$ are the Gabriel-Zisman homology $H_n(B\mC, \mG\circ \delta^{op})$, $n\geq 0$.
\end{definition}

Theorem \ref{homoliso} and Lemma \ref{confBW} lead to the following assertion about non-Abelian homology, which corresponds to the theorem about Abelian cohomology of the factorization category \cite[Theorem 4.4]{bau1985}.

\begin{corollary}\label{corfact}
For any contravariant natural system on an arbitrary small category $\mC$ there exist natural isomorphisms $H^{BW}_n(\mC, \mG)\cong \colim^{(\fF\mC)^{op}}_n\mG $, for all $n\geq 0$.
\end{corollary}

\subsection{Homotopy cofinality for non-Abelian 
Baues-Wirsching homology}

\begin{theorem}\label{confhomolBW}
Let $S: \mC\to \mD$ be a functor between small categories.
If for each $\alpha\in \Mor\mC$ the nerve of the category $S\langle\alpha\rangle$ is weakly equivalent to a point, then for every contravariant natural system $\mG$ on $\mD$ there are isomorphisms of non-Abelian Baues-Wirsching homology $H^{BW}_n(\mC, \mG\circ (\fF S)^{op}) \xrightarrow{\cong} H^{BW}_n(\mD, \mG)$
for all $n\geq 0$.
\end{theorem}
{\sc Proof.} Due to the isomorphism $\fF S \downarrow \alpha \cong \fF(S\langle \alpha\rangle)$ (Proposition \ref{factfibres}) and the contractibility of the nerve of the category $S\langle\alpha\rangle$, and the homotopy coinitiality of the functor $\fF S\langle\alpha\rangle \to S\langle\alpha\rangle$ (Proposition \ref{wefrac}), the category $\fF S\downarrow\alpha$ is contractible. Hence, we can apply Theorem \ref{homoliso} leading to isomorphisms $\colim^{(\fF\mC)^{op}}_n \mG(\fF S)^{op}
\xrightarrow{\cong} \colim^{(\fF\mD)^{op}}_n \mG$ and use Corollary \ref{corfact}.
\hfill$\Box$

\section{Conclusion}

The problems posed in the introduction are solved positively (Theorem \ref{homoliso} and Theorem \ref{discvirt}). When solving them, auxiliary statements for homotopy colimits of pointed simplicial sets are proved (Theorem \ref{cofpointed} and Proposition \ref{hoclan}).
The solution of these problems allowed us to construct the theory of non-Abelian homology of Gabriel-Zisman for simplicial sets (Corollary \ref{gzisos}), generalize it to $\mD$-sets with coefficients in group diagrams (Corollary \ref{dliso} and Corollary \ref{dhiso}), as well as to apply it to study non-Abelian Baues-Wirsching homology and functors between the categories of factorizations (Corollary \ref{corfact} and Theorem \ref{confhomolBW}).

\end{document}